\documentclass[11pt]{article}
\usepackage[fleqn]{amsmath}
\usepackage{amssymb}
\usepackage{graphicx}

\newtheorem{theorem}{Theorem}[section]
\newtheorem{corollary}{Corollary}[section]

\newtheorem{lemma}[theorem]{Lemma}

\def\x{$\hfill\rlap{$\sqcup$}\sqcap$\bigskip}

\def\Gal{\mathop{\rm Gal}\nolimits}

\begin{document}


\title{A Few More  Functions That Are Not APN Infinitely Often}

\author{ Yves Aubry\\
Institut de Math\'ematiques de Toulon\\
 Universit\'e du Sud Toulon-Var\\
  \bigskip
 France\\
 Gary McGuire\thanks{Research supported by the Claude
Shannon Institute, Science Foundation Ireland Grant 06/MI/006}\\
 School of Mathematical Sciences\\
 University College Dublin\\
 \bigskip
 Ireland\\
  Fran\c cois Rodier\\
  Institut de Math\'ematiques de Luminy\\
  C.N.R.S., Marseille\\
   France
  }


\date{}
 \maketitle

\begin{abstract}
\noindent 
We consider exceptional APN functions on 
$\mathbb{F}_{2^m}$, which by definition are functions that are APN
on infinitely many extensions of $\mathbb{F}_{2^m}$.
Our main result is that polynomial functions of odd degree are not exceptional,
provided the degree is not a Gold number ($2^k+1$) or a Kasami-Welch number ($4^k-2^k+1$).
We also have partial results on functions of even degree, and functions
that have degree $2^k+1$. 
\end{abstract}


\bigskip


\bigskip


\section{Introduction}

Let $L=\mathbb{F}_q$ with $q=2^n$ for some positive integer $n$.
A function $f : L \longrightarrow L$ is said to be \emph{almost perfect nonlinear} (APN) on $L$
if the number of solutions in $L$ of the equation
$$f(x+a)+f(x)=b$$
is at most 2, for all $a,b\in L$, $a \not=0$. Equivalently, $f$ is APN if the set
$\{f(x+a)+f(x): x \in L \}$ has size at least $2^{n-1}$ for each $a \in L^*$. 
Because $L$ has characteristic 2, the number of solutions
to the above equation must be an even number, for any function $f$ on $L$.

This kind of function is very useful in cryptography because of its good resistance to differential cryptanalysis as was  proved by  Nyberg in \cite{ny}.

The best known examples of APN functions are the Gold functions $x^{2^k+1}$ and
the Kasami-Welch functions $x^{4^k-2^k+1}$.
These functions are defined over $\mathbb{F}_{2}$, and are APN on 
any field $\mathbb{F}_{2^m}$ where $gcd(k,m)=1$.

If $f$ is APN on $L$, then $f$ is APN on any subfield of $L$ as well.
We will consider going in the opposite direction. 
Recall that every function $f : L \longrightarrow L$ can be expressed as a polynomial
with coefficients in $L$, and
this expression is unique if the degree is less than $q$.
We can ``extend" $f$ to an extension field of $L$ by
using the same unique polynomial formula to define a function
on the extension field.
With this understanding, we will consider functions $f$ which are APN on $L$,
and we ask whether $f$ can be APN on an extension field of $L$.
More specifically, we consider functions that are APN on 
infinitely many extensions of $L$.
We call a function $f : L \longrightarrow L$  \emph{exceptional} if $f$ is APN on $L$
and is also APN on infinitely many extension fields of $L$.
The Gold and Kasami-Welch functions are exceptional.

We make the following conjecture.

\bigskip

{\bf Conjecture:}  Up to equivalence, the Gold and Kasami-Welch  functions are the only 
exceptional APN functions.
\bigskip

Equivalence here refers to CCZ equivalence; for a definition and discussion of this
see \cite{bcp}  for example.

We will prove some cases of this conjecture. 
It was proved in Hernando-McGuire \cite{HM} that the conjecture is true among
the class of monomial functions.
Some cases for $f$ of small degree have been proved by Rodier \cite{FR}.

We define
\[
\phi (x,y,z) = \frac{f(x)+f(y)+f(z)+f(x+y+z)}{(x+y)(x+z)(y+z)}
\]
which is a polynomial in $\mathbb{F}_q [x,y,z]$.  This  polynomial defines a surface $X$ in the three dimensional affine space $ \mathbb{A}^3$.

If $X$ is absolutely irreducible (or has an absolutely irreducible component
defined over $\mathbb{F}_q$) then $f$ is not APN on 
$\mathbb{F}_{q^n}$ for all $n$ sufficiently large.
As shown in \cite{FR}, this follows from the Lang-Weil bound for surfaces, which guarantees many
$\mathbb{F}_{q^n}$-rational points on the surface for all $n$ sufficiently large.

Let $\overline X$ denote the projective closure of $X$ in the three dimensional projective space  $\mathbb{P}^3 $.
If $H$ is a another projective hypersurface in  $\mathbb{P}^3 $, the idea of this
paper is to apply the following lemma.

\begin{lemma}
 If  $\overline X \cap H$ is a reduced (no repeated component) absolutely irreducible curve, then 
 $\overline X$ is absolutely irreducible.
 \end{lemma}
 
 Proof:
If $\overline X$ is not absolutely irreducible then every irreducible component of $\overline X$ intersects $H$ in a variety of dimension at least 1 (see Shafarevich \cite[Chap. I, 6.2, Corollary 5]{sh}).
 So $\overline X \cap H$ is reduced or  reducible.
  
 \x
 
 In particular, we will apply this when $H$ is a hyperplane.
 In Section 2 we study functions whose degree is not a
 Gold number ($2^k+1$) or a Kasami-Welch number ($4^k-2^k+1$).
 In Section 3 we study functions whose degree is a Gold number - this case is more subtle.

 The equation of  $\overline X$ is the homogenization of $\phi(x,y,z)=0$, which
 is  $\overline \phi(x,y,z,t)=0$ say.
 If $f(x)=\sum_{j=0}^d a_j x^j$ write this as
 \[
 \overline \phi(x,y,z,t)=\sum_{j=3}^{d} a_{j} \phi_j (x,y,z) t^{d-j}
 \]
 where 
 \[
  \phi_j (x,y,z)=\frac{x^{j}+y^{j}+z^{j}+(x+y+z)^{j}}{(x+y)(x+z)(y+z)}
 \]
 is homogeneous of degree $j-3$.
 We will later consider the intersection of $\overline X$ with the hyperplane
 $z=0$, and this intersection is a curve in a two dimensional projective space   with equation
 $\overline \phi(x,y,0,t)=0$.  
 An affine equation of this surface $\overline X$ is $\overline \phi(x,y,z,1)=\phi (x,y,z)=0$.
 
 A fact we will use is that if $f(x)=x^{2^k+1}$ then
 \begin{equation}\label{Goldfactors}
 \phi(x,y,z)=\prod_{\alpha\in \mathbb{F}_{2^k}-\mathbb{F}_2}(x+\alpha y+(\alpha+1)z).
 \end{equation}
 This can be shown by elementary manipulations (see Janwa, Wilson, \cite[Theorem 4]{JW}).
 
 Our definition of exceptional APN functions is motivated by the definition
of exceptional permutation polynomials.
A permutation polynomial $f:\mathbb{F}_q \longrightarrow \mathbb{F}_q$
is said to be exceptional if $f$ is a permutation polynomial on infinitely
many extensions of $\mathbb{F}_q$.
One technique for proving that a polynomial is not exceptional is
to prove that the curve $\phi(x,y)=(f(y)-f(x))/(y-x)$
has an absolutely irreducible factor over $\mathbb{F}_q$.
Then the Weil bound applied to this factor guarantees many
$\mathbb{F}_{q^n}$-rational points on the curve for all $n$ sufficiently large.
In particular there are points with $x\not= y$, which means that $f$
cannot be a permutation.

The authors thank the referee for relevant suggestions.

 \section{Degree not Gold or Kasami-Welch}
 
 If the degree of $f$ is not a Gold number $2^k+1$, or a Kasami-Welch number
 $4^k-2^k+1$, then we will apply  results of Rodier \cite{FR} and Hernando-McGuire \cite{HM}
 to prove our results. 

\begin{lemma}\label{irred}
Let $H$ be a projective hypersurface.
If $\overline X \cap H$ has a reduced absolutely irreducible component defined over $\mathbb{F}_q$
then
$\overline X$ has an absolutely irreducible component defined over $\mathbb{F}_q$.
\end{lemma}

Proof:
Let $Y_H$ be a
 reduced absolutely irreducible component 
 of $\overline X \cap H$ defined over $\mathbb{F}_q$.
 Let $Y$ be an absolutely irreducible component of $\overline X$ that contains $Y_H$.
Suppose for the sake of contradiction that $Y$ is not defined over  $\mathbb{F}_q$.
Then $Y$ is defined over  $\mathbb{F}_{q^t}$ for some $t$.
Let $\sigma$ be a generator for the Galois group $\Gal(\mathbb{F}_{q^t}/\mathbb{F}_{q})$ of $\mathbb{F}_{q^t}$ over $\mathbb{F}_{q}$.
Then $\sigma(Y)$ is an absolutely irreducible component of $\overline X$ that is
distinct from $Y$.
However, $\sigma (Y) \supseteq \sigma (Y_H)=Y_H$, which implies that
$Y_H$ is contained in two distinct absolutely irreducible components of $\overline X$.
This means that a double copy of $Y_H$  is a component of $\overline X$, which contradicts the
assumption that $Y_H$ is reduced.
\x

\begin{lemma}\label{reduction}
Let $H$ be the hyperplane at infinity.
Let $d$ be the degree of $f$.
Then $\overline X\cap H$ is not reduced if $d$ is even, and $\overline X\cap H$ is reduced if $d$ is odd
and $f$ is not a Gold or Kasami-Welch monomial function.
\end{lemma}

Proof:
Let $\phi_d(x,y,z)$ denote the $\phi$ corresponding to the function $x^d$.
In  $\overline X\cap H$ we may assume $\phi=\phi_d$.

If $d$ is odd then the singularities of $\overline X\cap H$ were classified by Janwa-Wilson \cite{JW}.
They show that the singularities are isolated (the coordinates must be $(d-1)$-th
roots of unity) and so the dimension of the singular locus of $\overline X\cap H$ is 0.

Suppose $d$ is even and write $d=2^je$ where $e$ is odd.
In  $\overline X\cap H$ we have
\begin{eqnarray*}
(x+y)(x+z)(y+z)\phi_d(x,y,z)&=& x^d+y^d+z^d+(x+y+z)^d\\
&=& (x^e+y^e+z^e+(x+y+z)^e)^{2^j}\\
&=& ((x+y)(x+z)(y+z)\phi_e(x,y,z) )^{2^j}.
\end{eqnarray*}
Therefore
\[
\phi_d (x,y,z)=\phi_e(x,y,z)^{2^j} ((x+y)(x+z)(y+z))^{2^j-1}
\]
and is not reduced.
\x

Here is the main result of this section.

\begin{theorem}
If the degree of the polynomial function $f$ is odd and not a Gold or a Kasami-Welch number then $f$ is not APN over $\mathbb{F}_{q^n}$ for all $n$ sufficiently large. 
\end{theorem}

Proof:
By Lemma \ref{reduction}, $\overline X\cap H$ is reduced. Furthermore, we know by \cite{HM} that $\overline X\cap H$ has an absolutely irreducible component defined over $\mathbb{F}_q$, which is also reduced. Thus, by Lemma \ref{irred}, we obtain that $\overline X$ has an absolutely irreducible component defined over $\mathbb{F}_q$. As discussed in the introduction, this enables us to conclude that $f$ is not APN on 
$\mathbb{F}_{q^n}$ for all $n$ sufficiently large.
\x

In the even degree case, we can state the result when half of the degree is odd, with an extra minor condition.

\begin{theorem}
If the degree of the polynomial  function $f$ is $2e$
with $e$ odd,   and if $f$ contains a term of odd degree, then $f$ is not APN over $\mathbb{F}_{q^n}$ for all $n$ sufficiently large. 
\end{theorem}
Proof:
As shown in the proof of Lemma \ref{reduction} in the particular case where $d=2^je$ with $e$ odd and $j=1$, we can write
$$\phi_d (x,y,z)=\phi_e(x,y,z)^{2} (x+y)(x+z)(y+z).$$
Hence, $x+y=0$ is the equation of a reduced component of the curve $X_{\infty}=\overline X\cap H$ with equation $\phi_d=0$ where $H$ is the hyperplane at infinity. The only absolutely irreducible component $X_0$ of the surface $\overline X$ containing the line $x+y=0$ in $H$ is reduced and defined over $\mathbb{F}_q$. We have to show that this component doesn't contain the plane $x+y=0$.

The function $x+y$ doesn't divide $\phi(x,y,z)$ if and only if the function $(x+y)^2$ doesn't divide $f(x)+f(y)+f(z)+f(x+y+z)$.
Let $x^r$ be a term of odd degree of the function $f$. We show easily that   $(x+y)^2$ doesn't divide $x^r+y^r+z^r+(x+y+z)^r$ by using the change of variables $s=x+y$ which gives:
$$x^r+y^r+z^r+(x+y+z)^r=s(x^{r-1}+z^{r-1})+s^2P$$
where $P$ is a polynomial. 

Hence $\overline X$ has an absolutely irreducible component defined over $\mathbb{F}_q$ and then $f$ is not APN on 
$\mathbb{F}_{q^n}$ for all $n$ sufficiently large.

\x

\noindent
{\bf Remark}:
This theorem is false if $2e$ is replaced by $4e$ in the statement.
A counterexample is $x^{12}+cx^3$, where $c\in \mathbb{F}_{4}$ satisfies $c^2+c+1=0$,
which is APN on $\mathbb{F}_{4^n}$ for
any $n$ which is not divisible by 3, since it is CCZ-equivalent to $x^3$.
Indeed
this function is defined over $\mathbb{F}_{4}$,
and is equal to $L\circ f$, where $f(x)=x^3$ and $L(x)=x^4+cx$.
Certainly $L$ is $\mathbb{F}_{4}$-linear, and
it is not hard to show that $L$ is bijective on 
$\mathbb{F}_{4^n}$ if and only if 
$n$  is not divisible by 3.
The graph of $x^3$ is $\{(x,x^3)\mid x\in \mathbb{F}_{4^n}\}$ and it is transformed in the graph of $x^{12}+cx^3$ which is $\{(x,x^{12}+cx^3)\mid x\in \mathbb{F}_{4^n}\}$ by the linear permutation $Id\times L$ where $Id$ is the identity function.
So when $n$ is not divisible by 3, $L\circ f$ is APN on $\mathbb{F}_{4^n}$
because $f$ is APN. This example shows in particular that our conjecture has to be stated up to 
CCZ-equivalence.

\section{Gold Degree}

Suppose the degree of $f$ is a Gold number $d=2^k+1$. 
Set $d$ to be this value for this section.
Then the degree of $\phi$ is $d-3=2^k-2$.

\subsection{First Case}

We will prove the absolute irreducibility for a certain type of $f$.

\begin{theorem}
Suppose $f(x)=x^d+g(x)$ where $\deg (g) \leq 2^{k-1}+1$ .
Let $g(x)=\sum_{j=0}^{2^{k-1}+1} a_j x^j$.
Suppose moreover that there exists a nonzero coefficient $a_j$ of $g$ such that 
$\phi_{j} (x,y,z)$ is absolutely irreducible.
Then $\phi (x,y,z)$ is absolutely irreducible.
\end{theorem}

Proof:  
We must show that $\phi(x,y,z)$ is absolutely irreducible.
Suppose $\phi(x,y,z)=P(x,y,z)Q(x,y,z)$.
Write each polynomial as a sum of homogeneous parts:
\begin{equation}\label{homogeneous}
\sum_{j=3}^{d} a_{j}\phi_j (x,y,z)= (P_s+P_{s-1}+\cdots +P_0)(Q_t+Q_{t-1}+\cdots +Q_0)
\end{equation}
where $P_j, Q_j$ are homogeneous of degree $j$.
Then from (\ref{Goldfactors}) we get
\[
P_sQ_t=\prod_{\alpha\in \mathbb{F}_{2^k} - \mathbb{F}_2}(x+\alpha y+(\alpha+1)z).
\]
In particular this implies that  $P_s$ and $Q_t$ are relatively prime as the product is made of distinct irreducible factors.

The homogeneous terms in (\ref{homogeneous}) of degree strictly less than $d-3$ and strictly greater than $2^{k-1}-2$ are 0,
by the assumed bound on the degree of $g$.
Equating terms of degree $s+t-1$ in the equation (\ref{homogeneous}) gives
$P_s Q_{t-1}+P_{s-1}Q_t=0$.  Hence $P_s$ divides $P_{s-1}Q_t$
which implies $P_s$ divides $P_{s-1}$
because $gcd(P_s,Q_t)=1$, and we conclude $P_{s-1}=0$ as $\deg P_{s-1}<\deg P_{s}$.
Then we also get $Q_{t-1}=0$.
Similarly, $P_{s-2}=0=Q_{t-2}$, $P_{s-3}=0=Q_{t-3}$, and so on until we get the equation
\[
P_s Q_0+P_{s-t}Q_t=0
\]
where we suppose wlog that $s\geq t$.
(Note that when $s\geq t$, one gets from $s+t=d-3$ that $s\geq (d-3)/2$ and $t\leq (d-3)/2$, and the bound on $\deg (g)$ is chosen: $\deg(g)<t+3\leq 2^{k-1}+2$.)
This equation implies $P_s$ divides $P_{s-t}Q_t$, which implies $P_s$ divides
$P_{s-t}$, which implies $P_{s-t}=0$.
Since $P_s\not=0$ we must have $Q_0=0$.

We now have shown  that $Q=Q_t$ is  homogeneous.
In particular, this means that $\phi_j (x,y,z)$ is divisible by $x+\alpha y+(\alpha+1)z$ for
some $\alpha \in  \mathbb{F}_{2^k}- \mathbb{F}_2$ and for all $j$ such that $a_{j}\ne0$.
We are done if there exists such a $j$ with $\phi_j (x,y,z)$ irreducible.

\x

\noindent
{\bf Remark}:
The hypothesis that there should exist a $j$ with $\phi_{j} (x,y,z)$ is absolutely irreducible
is not a strong hypothesis.  This is true in many cases (see the next remarks).
However, some hypothesis is needed, because the theorem is false without it.
One counterexample is with $g(x)=x^5$ and $k\geq 4$ and even.

\noindent
{\bf Remark}:
It is known that $\phi_{j}$ is
irreducible in the following cases (see \cite{jmw}):
\begin{itemize}
\item  $j\equiv3 \pmod 4$;
\item  $j\equiv5 \pmod 8$ and $j > 13$.
\end{itemize}

\noindent
{\bf Remark}:
The theorem is true with the weaker hypothesis that 
 there exists a nonzero coefficient $a_j$ such that $\phi_{j} (x,y,z)$ is prime to $\phi_d$ (recall $d=2^k+1$).
 This is the case for
 \begin{itemize}
\item $j=2^r+1$ is a Gold exponent with $r$ prime to $k$;
\item $j$ is a Kasami exponent (see \cite[Theorem 5]{JW});
\item $j=2^je$ with $e$ odd and $e$ is in one of the previous cases.
\end{itemize}

\noindent
{\bf Example:}  This applies to 
$x^{33}+g(x)$ where $g(x)$ is any polynomial of degree $\leq 17$.

\noindent
{\bf Remark:}  The proof did not use the fact that $f$ is APN.  This is simply a result about polynomials.

\noindent
{\bf Remark:} The bound $\deg (g) \leq 2^{k-1}+1$ is best possible, in the sense
that there is an example with $\deg (g) = 2^{k-1}+2$ in Rodier  \cite{FR}
where $\phi$ is not absolutely irreducible.   
The counterexample has $k=3$, and
$f(x)=x^9+ax^6+a^2x^3$.
We discuss this in the next section.

\subsection{On the Boundary of the First Case}

As we said in the previous section, when $f(x)=x^{2^k+1}+g(x)$
with $\deg (g) = 2^{k-1}+2$, it is false that
$\phi$ is always absolutely irreducible.   
However, the polynomial $\phi$ corresponding to the counterexample
$f(x)=x^9+ax^6+a^2x^3$ where $a\in \mathbb{F}_q $
factors into two irreducible factors over $\mathbb{F}_q $.
We generalize this to the following theorem.

\begin{theorem}
Let $q=2^n$.
Suppose $f(x)=x^d+g(x)$ where $g(x)\in \mathbb{F}_q [x]$ and $\deg (g) =2^{k-1}+2$.
Let $k$ be odd and relatively prime to $n$.
If $g(x)$ does not have the form $ax^{2^{k-1}+2}+a^2x^3$ then
$\phi$ is absolutely irreducible, while if
$g(x)$ does have the form $ax^{2^{k-1}+2}+a^2x^3$
then either $\phi$ is irreducible or $\phi$ splits into two absolutely irreducible factors which are both
defined over $\mathbb{F}_q$.
\end{theorem}

Proof:  
Suppose $\phi(x,y,z)=P(x,y,z)Q(x,y,z)$ and
let 
$$g(x)=\sum_{j=0}^{2^{k-1}+2} a_j x^j.$$
Write each polynomial as a sum of homogeneous parts:
\[
\sum_{j=3}^{d} a_{j}\phi_j (x,y,z)= (P_s+P_{s-1}+\cdots +P_0)(Q_t+Q_{t-1}+\cdots +Q_0).
\]
Then
\[
P_sQ_t=\prod_{\alpha\in \mathbb{F}_{2^k}-\mathbb{F}_2}(x+\alpha y+(1+\alpha )z).
\]
In particular this means $P_s$ and $Q_t$ are relatively prime as in the previous theorem.
We suppose wlog that $s\geq t$, which implies $s\geq 2^{k-1}-1$.
Comparing each degree gives
$P_{s-1}=0=Q_{t-1}$, $P_{s-2}=0=Q_{t-2}$, and so on until we get the equation
of degree $s+1$
\[
P_{s} Q_1+P_{s-t+1}Q_{t}=0
\]
which implies $P_{s-t+1}=0=Q_1$.
If $s\not=t$ then $s\geq 2^{k-1}$.
Note then that $a_{s+3} \phi_{s+3}=0$.
The equation of degree $s$ is
\[
P_{s} Q_0+P_{s-t}Q_{t}=a_{s+3} \phi_{s+3}=0.
\]
This means that $P_{s-t}=0$, so $Q_0=0$.
We now have shown  that $Q=Q_t$ is  homogeneous.
In particular, this means that $\phi(x,y,z)$ is divisible by $x+\alpha y+(1+\alpha )z$ for
some $\alpha \in  \mathbb{F}_{2^k}- \mathbb{F}_2$, which is impossible.
Indeed, since the leading coefficient of $g$ is not 0, the polynomial $\phi_{2^{k-1}+2}$ occurs in $\phi$;
as 
$\phi_{2^{k-1}+2}= \phi^2_{2^{k-2}+1}(x+y)(y+z)(z+x)$, this polynomial is prime to $\phi$, because if $x+\alpha y+ (1+\alpha)z$ occurs  in the two polynomials $\phi_{2^{k-1}+2}$ and $\phi_{2^k+1}$, then $\alpha$ would be an element of  $\mathbb{F}_{2^k}\cap\mathbb{F}_{2^{k-2}}=\mathbb{F}_{2}$ because $k$ is odd.

Suppose next that
$s=t=2^{k-1}-1$ in which case
the degree $s$ equation is
\[
P_s Q_0+P_{0}Q_s=a_{s+3} \phi_{s+3}.
\]

If $Q_0=0$, then
$$\phi(x,y,z)=\sum_{j=3}^{d} a_{j}\phi_j (x,y,z)= (P_s+P_0) Q_t$$
which implies that
$$\phi(x,y,z)=a_{d}\phi_d(x,y,z)+a_{2^{k-1}+2} \phi_{2^{k-1}+2} (x,y,z)= P_sQ_t+P_0 Q_t$$
and  $P_0\ne0$, since $g\ne0$.
So one has
$\phi_{2^{k-1}+2} $ divides $\phi_d(x,y,z)$
which is impossible as
$$\phi_{2^{k-1}+2} =\phi_{2^{k-2}+1} ^2 (x+y)(y+z)(z+x).$$

We may assume then that $P_0=Q_0$, and we have $\phi_{2^{k-1}+2} =0$.
Then we have
\begin{equation}\label{phifactors}
\phi(x,y,z)=(P_s+P_0)(Q_s+Q_0)=P_sQ_s+P_0(P_s+Q_s)+P_0^2.
\end{equation}
Note that this implies $a_j=0$ for all $j$ except $j=3$ and $j=s+3$.  This means
\[
f(x)=x^d+a_{s+3}x^{s+3}+a_3x^3.
\]
So if $f(x)$ does not have this form, this shows that $\phi$ is absolutely irreducible.

If on the contrary $\phi$ splits as $(P_s+P_0)(Q_s+Q_0)$, the factors $P_s+P_0$ and $Q_s+Q_0$ are irreducible, as can be  shown by using the same argument.

Assume from now on that $f(x)=x^d+a_{s+3}x^{s+3}+a_3x^3$ and
that (\ref{phifactors}) holds.
Then $a_3=P_0^2$,
so clearly $P_0=\sqrt{a_3}$ is defined over $\mathbb{F}_q$.
We claim that $P_s$ and  $Q_s$  are actually defined over $\mathbb{F}_2$.


We know from (\ref{Goldfactors}) that $P_sQ_s$ is defined over $\mathbb{F}_2$.

Also $P_0(P_s+Q_s)=a_{s+3} \phi_{s+3}$, 
so $P_s+Q_s=(a_{s+3}/\sqrt{a_3}) \phi_{s+3}$.
On the one hand, $P_s+Q_s$ is defined over $\mathbb{F}_{2^k}$ by (\ref{Goldfactors}). 
On the other hand, 
since $\phi_{s+3}$ is defined over $\mathbb{F}_2$ 
we may say that $P_s+Q_s$ is defined over $\mathbb{F}_q$.
Because $(k,n)=1$ we may conclude that
$P_s+Q_s$ is defined over $\mathbb{F}_2$.
Note that the leading coefficient of $P_s+Q_s$ is 1, so $a_{s+3}^2=a_3$.
Whence if this condition is not true, then $\phi$ is absolutely irreducible.

Let $\sigma$ denote the Galois automorphism $x\mapsto x^2$.
Then $P_sQ_s=\sigma(P_sQ_s)=\sigma(P_s)\sigma(Q_s)$,
and $P_s+Q_s=\sigma(P_s+Q_s)=\sigma(P_s)+\sigma(Q_s)$.
This means $\sigma$ either fixes both $P_s$ and $Q_s$, in which case we are done,
or else $\sigma$ interchanges them.
In the latter case, $\sigma^2$ fixes both $P_s$ and $Q_s$, so they are defined over
$\mathbb{F}_4$.
Because they are certainly defined over $\mathbb{F}_{2^k}$ by (\ref{Goldfactors}),
and $k$ is odd, they are defined over 
$\mathbb{F}_{2^k} \cap \mathbb{F}_{4} =\mathbb{F}_{2}$. 

Finally, we have now shown that $\overline X $ either is irreducible, or splits into two absolutely irreducible
factors defined over $\mathbb{F}_q$.
\x

\subsection{Using the Hyperplane $y=z$}

We study the intersection of  $\phi(x,y,z)=0$ with the hyperplane $y=z$.

\begin{lemma}
$\phi(x,y,y)$ is always a square.
\end{lemma}

Proof:
It suffices to prove the result for $f(x)=x^d$.
This is equivalent to proving that $\phi_d(x,1,1)$ is a square.
This is equivalent to showing that its derivative with respect to $x$ is identically 0.
This is again equivalent to showing that the partial derivative with respect to $x$
of $\phi_d (x,y,1)$, evaluated at $y=1$, is 0.
In Lemma 4.1 of \cite{FR} Rodier proves that $y+z$ divides the partial derivative
of $\phi_d (x,y,z)$ with respect to $x$, which is exactly what is required.
\x

\begin{lemma}
Let $H$ be the hyperplane $y=z$.
If $\overline X \cap H$ is the square of an absolutely irreducible component 
defined over $\mathbb{F}_q$ then
$\overline X$ is absolutely irreducible.
\end{lemma}

Proof: 
We claim that for any nonsingular point $P\in \overline X \cap H$, 
the tangent plane to the curve $\overline X \cap H$ at $P$ is $H$.
The equation of the tangent plane is 
\[
(x-x_0) \phi'_x (P) +(y-y_0) \phi'_y (P)+(z-z_0) \phi'_z (P)=0
\]
where $P=(x_0,y_0,z_0)$.
Since $P\in H$ we have $y_0=z_0$.
It is straightforward to show that $\phi'_x (P)=0$ and 
$\phi'_y (P)=\phi'_z (P)$, so this equation becomes
\[
(y+z) \phi'_y (P)=0.
\]
But $y+z=0$ is the equation of $H$.
\x

\begin{corollary} If $f(x)=x^d +g(x)$, and $d=2^k+1$ is a Gold exponent,
and $\phi(x,y,y)$ is the square of an irreducible, then $\overline X$ is absolutely irreducible.
\end{corollary}

Note that any term $x^d$ in $g(x)$ where $d$ is even will drop out when
we calculate $\phi(x,y,y)$, because if $d=2e$ then
\begin{eqnarray*}
\phi_d (x,y,z)&=& \frac{x^d+y^d+z^d+(x+y+z)^d}{(x+y)(x+z)(y+z)}\\
&=&\frac{(x^e+y^e+z^e+(x+y+z)^e)^{2}}{(x+y)(x+z)(y+z)}\\
&=& \phi_e (x,y,z) (x^e+y^e+z^e+(x+y+z)^e)\\
&=&0 \quad \textrm{on $H$}
\end{eqnarray*}
because the right factor vanishes on $H$.

In order to find examples of where we can apply this Corollary, if we write
\[
\phi(x,y,y)=(x+y)^{2^k-2}+h(x,y)^2
\]
then to apply this result we want to show that
\[
(x+y)^{2^{k-1}-1} + h(x,y)
\]
is irreducible.  The degree of $h$ is smaller than $2^{k-1}-1$.
Letting $t=x+y$ we want an example of $h$ with $t^{2^{k-1}-1}+h(x,x+t)$ is irreducible.

{\bf Example:} Choose any $h$ so that
 $h(x,x+t)$ is a monomial, and then  $t^{2^{k-1}-1}+h(x,x+t)$ is irreducible.


\end{document}